\documentclass[11pt]{article}
\usepackage{amssymb}
\title{Probabilistic characterisation of Besov-Lipschitz spaces on  metric measure spaces}
\author{Katarzyna Pietruska-Pa\l{}uba\thanks{Supported by a KBN grant no.
                         1-PO3A-008-29. }\\
                         Institute of Mathematics\\
                         University of Warsaw\\
                         ul. Banacha 2\\
                         02-097 Warsaw, Poland\\
                         e-mail: kpp@mimuw.edu.pl}
\renewcommand{\em}{\sl}

plus 6pt minus 12pt
\newcommand{\barint}{
         \rule[.036in]{.12in}{.009in}\kern-.16in
          \displaystyle\int  }
\def\r{{\mathbb{R}}}

\begin{document}

\newtheorem{theo}{\bf Theorem}[section]
\newtheorem{coro}{\bf Corollary}[section]
\newtheorem{lem}{\bf Lemma}[section]
\newtheorem{rem}{\bf Remark}[section]
\newtheorem{defi}{\bf Definition}[section]
\newtheorem{ex}{\bf Example}[section]
\newtheorem{fact}{\bf Fact}[section]
\newtheorem{prop}{\bf Proposition}[section]

\makeatletter \@addtoreset{equation}{section}
\renewcommand{\theequation}{\thesection.\arabic{equation}}
\makeatother

\newcommand{\ds}{\displaystyle}
\newcommand{\ts}{\textstyle}
\newcommand{\ol}{\overline}
\newcommand{\wt}{\widetilde}
\newcommand{\ck}{{\cal K}}
\newcommand{\ve}{\varepsilon}
\newcommand{\vp}{\varphi}
\newcommand{\pa}{\partial}
\newcommand{\rp}{\mathbb{R}_+}

\date{}
\maketitle

\begin{abstract}
\noindent We give a probabilistic characterisation of the
Besov-Lipschitz spaces $Lip(\alpha,p,q)(X)$ on domains which
support a Markovian kernel with appropriate exponential bounds.
This extends former results of \cite{Jon,KPP1,KPP2,GHL} which were
valid for $\alpha=\frac{d_w}{2},p=2$, $q=\infty,$ where $d_w$ is
the walk dimension of the space $X.$
\end{abstract}
\medskip

\section{Introduction}

There are several definitions of Besov-Lipschitz spaces on measure
spaces.  In this paper we will investigate the spaces
$Lip(\alpha,p,q)(X),$ as defined in Jonsson \cite{Jon}. Jonsson's
paper dealt with the Sierpi\'nski gasket embedded in
$\mathbb{R}^d$ only, but did not really use the embedding itself,
and so this particular definition can be extended to general
metric measure spaces (see e.g. \cite{Gri}, \cite{Kum}). The most
convenient to analyse are those spaces  on which there exists a
complete, symmetric Markovian kernel with appropriate exponential
bounds. There are several results concerning such spaces, see e.g.
\cite{GHL}, \cite{KPP2}, \cite{KPP3}, \cite{HZ}.

The existence of a Markovian kernel on $X$ of this type is
equivalent to the existence of a {\em fractional diffusion} on X
(see \cite{Bar} for the definition). Its generator, often called
`the Laplacian' on a general metric space, serves as a substitute
for the {\em bona fide} differentiation operator, even though the
differential itself is not-so-convenient to define in this
generality. This is one of the reasons why the existence of such a
kernel allows to prove certain properties of underlying spaces. In
particular, in a series of papers(\cite{Jon,KPP1,KPP2,GHL}) it has
been proven that the spaces $Lip(\frac{d_w}{2},2,\infty)(X)$
(where $d_w$ is the {\em walk dimension} of $X$) are domains of
the Dirichlet form associated with this particular diffusion --
and so it can be described using the kernel
$p(\cdot,\cdot,\cdot).$ Also, the spaces $Lip(\alpha,2,2)$ allow
for a probabilistic characterisation (see \cite{Sto}).
 In present work we extend these results and provide a
characterisation of the spaces $Lip(\alpha,p,q)(X)$ in terms of
the Markovian kernel whose existence we are assuming, for
$\alpha>0,$ $p,q\geq 1$ (nox excluding $q=\infty$). Our proof is
entirely elementary and uses only a variant of discrete Hardy
inequality (proven below).

Besov spaces, on the very same class of metric measure spaces,
were also introduced bu Hu and Z\"{a}hle -- in a different way -- in
their paper \cite{HZ}.  The way they are defined owes to the
classical characterisation of Besov spaces from \cite{Fle},
\cite{ST}. Those spaces will be denoted by $B^{p,q}_{\beta}(X);$
We will see that our characterisation of $Lip(\alpha,p,q)(X)$ is
consistent with $B^{p,q}_\beta(X)$ for some range of parameters
(see Section \ref{huza}).

In the case of simple fractals (and the Sierpi\'nski gasket in
particular), yet another definition of Besov-Lipschitz spaces was
given by Strichartz in \cite{Str}. This definition uses a discrete
approximation of the space $X.$ The Strichartz spaces we think of
are the spaces $(\Lambda^{p,q}_\alpha)^1(X)$ (in \cite{Str}, one
can find other spaces as well, corresponding to large values of
$\alpha$). Strichartz spaces do not allow the smoothness parameter
$\alpha$ to be too small; the definition uses a discrete
approximation of simple fractals, and so it is mandatory that the
functions concerned be continuous. This is not necessarily true
for small values of $\alpha$. It is known  (see \cite{Bod}) that
the Strichartz spaces and the Jonsson spaces agree for certain
range of parameters. Therefore our characterisation remains valid
for Strichartz spaces as well (see Section \ref{stri}).

\section{Preliminaries}

\noindent {\small{\bf Convention}. In the sequel, $c$ will denote
a generic constant whose value is irrelevant and can change from
line to line. The important constants will be denoted by upper
case letters or by letters with subscripts: $c_1,2,...$. When we
write $A\asymp B,$ then we mean that for some $c>0,$ $c^{-1}A\leq
B\leq cA.$}

 Suppose $(X,\rho)$ is a
locally compact metric space and that $\mu$ is a Borel measure on
$X$ which is Ahlfors $d-$regular, i.e. such that
\begin{equation}\label{measure}
\forall_{x\in X}\forall_{0<r<\mbox{\small diam}\,X}\;\;\;\;\;\;C_1
r^d\leq \mu(B(x,r))\leq C_2 r^d
\end{equation}
where $d>0$ and $C_1, C_2$ are positive constants.

Following Jonsson, we define the Besov-Lipschitz spaces as
follows.

\begin{defi} Suppose $\alpha>0,$ $p,q\in[1,\infty).$ Then $Lip(\alpha, p,q)(X)$
is, by definition, the collection of those $f\in L^p(X,\mu)$ for
which $\|(a_m(f))\|_q <+\infty,$ where
\begin{equation}\label{amf}
a_m(f)=2^{m\alpha}\left(2^{md}\int\int_{\rho(x,y)\leq
2^{-m}}|f(x)-f(y)|^p\,d\mu(x)d\mu(y)\right)^{1/p}.
\end{equation}
\end{defi}

 When $p=\infty$ or $q=\infty$ then the usual modifications are
 needed.

 The expression
 \[\|f\|_{\alpha,p,q}=\|(a_m(f))\|_q\]
is a seminorm, which can be turned into a norm by adding
$\|f\|_{L^p}.$ The space $Lip(\alpha,p,q)(X),$ equipped with the
norm
\begin{equation}\label{norm1}
\|f\|_{Lip(\alpha,p,q)}= \|f\|_{L^p}+ \|f\|_{\alpha,p,q}
\end{equation}
is a Banach space.

In the sequel, we  assume  that there exists  a symmetric
Markovian kernel   $\{p(t,x,y)\}_{t>0}$ on $X,$ i.e. a family of
measurable functions $p(t,\cdot,\cdot):X\times X\to \r_+,$ which
satisfies:
\begin{description}
  \item[(A1)] $\forall_{t>0}\forall_{x,y\in X}\;\;p(t,x,y)=p(t,y,x)$ (symmetry),
  \item[(A2)]  $\forall_{t>0}\forall_{x\in X}\;\;\int_X p(t,x,y) d\mu(y) =1$ (normalisation
  or stochastic completeness),
  \item[(A3)] $\forall_{s,t>0}\forall_{x,y\in X}\;\;p(s+t,x,y)=\int_X p(s,x,z)p(t,z,y) d\mu(z)$ (the
  Chapman-Kolmogorov identity, or the Markov property),
  \item[(A4)] $\forall_{t>0} \forall_{x,y\in X}\;\; p(t,x,y)>0$
  (irreducibility),
  \item[(A5)] $\forall f\in L^2(X),\;\; P_t f\to f$ when $t\to 0,$
  strongly in $L^2(X),$ where $P_tf(x)=\int_X f(y)
  p(t,x,y)\,d\mu(y)$ (strong continuity).
\end{description}

These conditions allow us to freely use the Dirichlet form theory
for Markov processes.

And, finally, our main assumption:
\begin{description}\item[(A6)]
\begin{equation}\label{esti2}
\frac{c_1}{t^{d/d_w}}e^{-c_2\left(\frac{\rho(x,y)}{t^{1/d_w}}\right)^{\frac{d_w}{d_w-1}}}
\leq p(t,x,y) \leq
\frac{c_3}{t^{d/d_w}}e^{-c_4\left(\frac{\rho(x,y)}{t^{1/d_w}}\right)^{\frac{d_w}{d_w-1}}}.
\end{equation}
\end{description}
The parameter $d_w$ is usually called the {\em walk dimension} of
$X,$ as it
 controls the weak time/space  scaling of the Markov process with
 transition density $p(\cdot,\cdot,\cdot)$.
 Such a process will be denoted by $(B_t,\mathbb{P}_x)_{t\geq 0, x\in X,}.$
We know that the parameter $d_w$ is the same for all possible
Markov  processes sharing the estimate (\ref{esti2}). When the
space $(X,\rho)$ satisfies the {\em chain condition}:

\begin{description}
\item[(CC)] there exists a constant $C>0$ such that for any ${x,y\in
X}$ and any positive integer $n$, there
  exists  a chain $x=x_0,x_1,...,x_n=y$
 of points from $X$
s.t. $\rho(x_i,x_{i+1})\leq\frac{C}{n}\rho(x,y),$
\end{description}
then
\[d\leq d_w\leq d+1.\]

  The estimate (\ref{esti2}) ensures that this process is a
  diffusion.
  In fact it is known (see \cite{Gri-Kum}) that when the process
  is a diffusion, then the
  only possible function $\Phi$ in the estimate of the form
  $\frac{c}{t^{d/d_w}}\Phi\left(\frac{\rho(x,y)}{t^{{1/d_w}}}\right)$
 is the exponential function as in (\ref{esti2}).

Among the examples, we can list
\begin{description}
\item[$\bullet$]  the Euclidean space $\mathbb{R}^d$ with the Gaussian kernel, $g(t,x,y)= \frac{1}{(2\pi
t)^{d/2}} \exp(-\frac {|x-y|^2}{2t}),$ certain manifolds with
nonnegative Ricci curvature,
\item[$\bullet$] simple fractals, where  $p(t,x,y)$
is the transition density of the Brownian motion, and can be
bounded from both above and below by
$\frac{c}{t^{d/d_w}}\exp(-c(\frac{\rho(x,y)}{t^{1/d_w}})^{d_w/(d_w-1)}),$
see \cite{Kum0},
\item[$\bullet$] the Brownian motion on p.c.f. self similar  sets and
on the Sierpi\' nski carpets, where we have an estimate analogous
to that on simple nested fractals (see \cite{Ham-Kum} and
\cite{Bar-Bas}).
\end{description}

\section{The main theorem}
In a series of papers (\cite{Jon,KPP1,KPP2,GHL}) it has been
proven that the domain of the Dirichlet form associated with the
Markovian kernel satisfying {\bf (1)-(6)} is equal to the space
$Lip(\frac{d_w}{2}, 2,\infty)(X),$ and that the norms: of the
Dirichlet space ${\cal D}({\cal E})$ and of the Besov-Lipschitz
space $Lip(\frac{d_w}{2},2,\infty)(X)$ are equivalent.

The domain of the Dirichlet form, ${\cal D}({\cal E}),$ consists
of those functions $f\in L^2(X,\mu),$ for which $s(f)<\infty,$
where
\begin{equation}\label{sf}
s(f)=\sup_{t>0} \frac{1}{2t}\int_X\int_X
(f(x)-f(y))^2p(t,x,y)\,d\mu(x)d\mu(y). \end{equation}

 We will give
a similar characterisation of the spaces $Lip(\alpha,p,q)(X),$ for
general parameters $\alpha >0,$ $p,q\geq 1.$ Namely, we show:

\begin{theo}\label{mmain}
Suppose  $\alpha >0$ and $ p,q\in[1,\infty).$  Then $f$ belongs to
the Lipschitz-Besov space $Lip(\alpha,p,q)$ if and only if:\\
{\bf (1)} $f\in L^p(X,\mu),$\\
{\bf (2)} \begin{eqnarray*}
  I^{(\alpha)}(f)&:= & \int_0^1 \frac{1}{t^{\frac{\alpha
  q}{d_w}}} \left(\int_X\int_X |f(x)-f(y)|^p p(t,x,y) \,d\mu(x)d\mu(y)
  \right)^{\frac{q}{p}}\frac{dt}{t}<\infty.
  \end{eqnarray*}
  Moreover, we have
  \begin{equation}\label{compare}
  \|f\|_{Lip(\alpha,p,q)}\asymp
  \|f\|_{L^p}+\left(I^{(\alpha)}(f)\right)^{\frac{1}{q}}.
  \end{equation}
\end{theo}

\noindent {\bf Proof.} \textsc{Part 1.} Suppose that $f\in
L^p(X,\mu)$ and that $I^{(\alpha)}(f)<\infty.$

For later use, introduce the notation
\begin{equation}\label{imale}
i_m^{(\alpha)}(f)=\int\int_{\rho(x,y)\leq
2^{-m}}|f(x)-f(y)|^p\,d\mu(x)d\mu(y).
\end{equation}

Clearly, we have
\[ I^{(\alpha)}(f)\geq \int_0^1 \frac{1}{t^{\frac{\alpha q}{d_w}}}
\left(\int\int_{\rho(x,y)\leq t^{\frac{1}{d_w}}}|f(x)-f(y)|^p
p(t,x,y)\,d\mu(x)d\mu(y)\right)^{\frac{q}{p}}\,\frac{dt}{t}.
\]
When $\rho(x,y)\leq t^{\frac{1}{d_w}},$ then $p(t,x,y)$ is nearly
constant, and so (\ref{esti2}) gives

\begin{eqnarray*}
I^{(\alpha)}(f)&\geq & c\int_0^1 \frac{1}{t^{\frac{\alpha
q}{d_w}}}\left(\frac{1}{t^{\frac{d}{d_w}}}\int\int_{\rho(x,y)\leq
t^{\frac{1}{d_w}}} |f(x)-f(y)|^p
\,d\mu(x)d\mu(y)\right)^{\frac{q}{p}}\frac{dt}{t}\\[2mm]
& = & c\sum_{m=0}^\infty \int_{2^{-(m+1)d_w}}^{2^{-md_w}}
\frac{1}{t^{\frac{\alpha
q}{d_w}}}\left(\frac{1}{t^{\frac{d}{d_w}}}\int\int_{\rho(x,y)\leq
t^{\frac{1}{d_w}}} |f(x)-f(y)|^p
\,d\mu(x)d\mu(y)\right)^{\frac{q}{p}}\frac{dt}{t}\\[2mm]
&\geq & c\sum_{m=0}^\infty 2^{m\alpha q}
(2^{md}i_m^{(\alpha)}(f))^{\frac{q}{p}} =
c\,\|(a_m^{(\alpha)}(f))\|_{\ell^q}^q = c\|f\|_{\alpha,p,q}.
\end{eqnarray*}
This proves the inequality
\[\|f\|_{L^p}+(I^{(\alpha)}(f))^{1/q}\geq c \|f\|_{Lip(\alpha,p,q)}.\]

\noindent\textsc{Part 2.} Now we prove the opposite inequality.
Suppose that $f\in Lip(\alpha,p,q)(X).$ Similarly as before, write
 \[
I^{(\alpha)}(f)=\sum_{m=0}^\infty \int_{2^{-(m+1)d_w}}^{2^{-mdw}}
\frac{1}{t^{\frac{\alpha q}{d_w}}} \left( \int_X\int_X
|f(x)-f(y)|^p p(t,x,y)\,d\mu(x)d\mu(y)
\right)^{\frac{q}{p}}\,\frac{dt}{t}.
 \]

For given $m,$ we split the inner
integral into two parts: over the
set $Z_m=\{\rho(x,y)\leq 2^{-m/2}\},$ and over its complement
$Z_m^c=\{\rho(x,y)> 2^{-m/2}\},$ i.e.

\begin{eqnarray*}
I^{(\alpha)}(f)&\leq& c\sum_{m=0}^\infty
\int_{2^{-(m+1)d_w}}^{2^{-md_w}} \frac{1}{t^{\frac {\alpha
q}{d_w}}} \left(\int\int_{Z_m} |f(x)-f(y)|^p
p(t,x,y)\,d\mu(x)d\mu(y)\right)^{\frac{q}{p}}\frac{dt}{t}\\[3mm]
&& +c\sum_{m=0}^\infty \int_{2^{-(m+1)d_w}}^{2^{-md_w}}
\frac{1}{t^{\frac {\alpha q}{d_w}}} \left(\int\int_{Z_m^c}
|f(x)-f(y)|^p
p(t,x,y)\,d\mu(x)d\mu(y)\right)^{\frac{q}{p}}\frac{dt}{t},\nonumber \\[3mm]
&=:& I_1^{(\alpha)}+I_2^{(\alpha)}.\nonumber
\end{eqnarray*}
 The integral over $Z_m^c$ is not bigger than (use symmetry)
\begin{eqnarray*}
2^p\int\int_{Z_m^c} |f(x)|^p p(t,x,y)\,d\mu(x)d\mu(y) & = & 2^p
\int_X |f(x)|^p \left( \int_{\{y:\rho(x,y)> 2^{-m/2}\}} p(t,x,y)
\,d\mu(y)\right)d\mu(x)\\[2mm]
&=& 2^p \|f\|_{L^p}^p \sup_{x\in X}\mathbb{P}_x[\rho(x, B_t)>
2^{-m/2}]\\
&\leq & 2^p \|f\|_{L^p}^p \exp \left(-c2^{(-m/2)
t^{-1/d_w}}\right)^{\frac{d_w}{d_w-1}}
\end{eqnarray*}
(we have used the estimate $\mathbb{P}_x[\rho(x,B_t)\geq
\delta]\leq \exp(-c(\delta t^{-1/d_w})^{\frac{d_w}{d_w-1}})$,
valid under our assumption {\bf (A6)}, see \cite{Bar}).

Observe that while integrating over $Z_m^c,$ the values of $t$ are
confined to $[2^{-(m+1)d_w}, 2^{-md_w}],$ and so the integral we
are estimating does not exceed
\[c\|f\|_{L^p}^p \exp(-c(2^{m/2})^{\frac{d_w}{d_w-1}}),\]
thus

\begin{eqnarray*}
I_2^{(\alpha)}(f) & \leq & c\sum_{m=0}^\infty
\int_{2^{-(m+1)d_w}}^{2^{-m}d_w}\frac{1}{t^{\frac{\alpha
q}{d_w}}}\left( \|f\|_{L_p}^p \exp
(-c(2^{m/2})^{\frac{d_w}{d_w-1}})
\right)^{\frac{q}{p}}\frac{dt}{t}\\
&\leq & c\|f\|_{L_p}^q \sum_{m=0}^\infty
2^{mq\alpha}\exp(-c(2^{m/2})^{\frac{d_w}{d_w-1}})\leq
c\|f\|_{L^p}^q.
 \end{eqnarray*}

 We are left with estimating $I_1^{(\alpha)}(f),$ which requires
 subtler tools.

The double integral over the set $Z_m$ can be written as
\begin{equation}\label{num}
\sum_{k=m/2}^\infty \int\int_{2^{-(k+1)}<\rho(x,y)\leq 2^k}
|f(x)-f(y)|^p p(t,x,y) d\mu(x)d\mu(y),
\end{equation}
and again,  in this integral we have $t\in [2^{-(m+1)d_w},
2^{-md_w}],$ so from the basic estimate (\ref{esti2}) for the
transition density, we get that (\ref{num}) is not bigger than
\begin{eqnarray*}
&&c\,\sum_{k=m/2}^\infty 2^{md}\left(\int\!\!\int_
{2^{-(k+1)}<\rho(x,y)\leq 2^{-k}} |f(x)-f(y)|^pd\mu(x)d\mu(y)
\right)\exp(-c_4 (2^m
2^{-k})^{\frac{d_w}{d_w-1}})\\
&\leq & c\,\sum_{k=m/2}^\infty 2^{md} \exp(-c_4 (2^{m-k}
)^{\frac{d_w}{d_w-1}})\, i_k^{(\alpha)}(f)
\end{eqnarray*}
($i_k^{(\alpha)}(f)$ was defined by (\ref{imale})). Consequently,
\begin{eqnarray}
I_1^{\alpha}(f)&\leq & c\,\sum_{m=0}^\infty 2^{m\alpha q} \left(
\sum_{k=m/2}^\infty
2^{md}\, \exp({-c_4(2^{(m-k)\frac{d_w}{d_w-1}})})\,i_k^{(\alpha)}(f) \right)^{\frac{q}{p}}\nonumber\\
&=&c\,\sum_{m=0}^\infty 2^{m q(\alpha+\frac{d}{p})}
 \left(
\sum_{k=m/2}^\infty i_k^{(\alpha)}(f)\,
\exp(-c_4(2^{(m-k)\frac{d_w}{d_w-1}})) \right)^{\frac{q}{p}}\nonumber\\
&\leq & c\,\sum_{m=0}^\infty 2^{m q(\alpha+\frac{d}{p})}
\left(\sum_{k=m/2}^m \dots\right)^{\frac{q}{p}} +
c\,\sum_{m=0}^\infty 2^{m q(\alpha+\frac{d}{p})}
\left(\sum_{k=m}^\infty \dots\right)^{\frac{q}{p}}\label{ttt}.
\end{eqnarray}
To estimate these double sums we will use the discrete Hardy
inequalities:  the classical Hardy inequality (\ref{clashardy})
for the second sum, and the modified Hardy inequality
(\ref{modhardy}) for the first one.

\noindent We include them as lemmas. Lemma \ref{class} is
classical so we omit its proof.

\begin{lem}[classical discrete Hardy inequality]\label{class}
Suppose $r>0.$ $t>1,$ $x_m>0,$ $m=1,2,...$ Then
\begin{equation}\label{clashardy}
\sum_{m=0}^\infty t^m\left(\sum_{k=m}^\infty x_m\right)^r\leq
K\,\sum_{m=0}^\infty t^m x_m^r,
\end{equation}
where  $K$  is a  constant depending on $r,t$ only.
\end{lem}

\begin{lem}[modified discrete Hardy inequality]\label{lemhardy}
Suppose $r>0,$ $t>1,$ $x_m>0,$ $\kappa >1,$ $\lambda>0.$ Then
\begin{equation}\label{modhardy}
\sum_{m=0}^\infty t^m \left(\sum_{k=m/2}^m
x_ke^{-\lambda\kappa^{m-k}}\right)^r\leq K \,\sum_{m=0}^\infty
 t^m x_m^r,
 \end{equation}
 where the constant $K$ depends  on $r,t,\kappa,\lambda$ only.
 \end{lem}

\noindent {\bf Proof of Lemma \ref{lemhardy}.}

\noindent \textsc{{Case 1. $r\leq 1.$}} Starting with the
elementary inequality
\[ (y_1+...+y_n)^r\leq y_1^r+...+y_n^r,\]
which is valid for all $n=1,2,...,$ $y_1,...,y_n \geq 0$ and
$r\in(0,1],$ we have
 \begin{eqnarray*}
 \sum_{m=0}^\infty t^m \left(\sum_{k=m/2}^m x_k
 \exp(-\lambda \kappa^{m-k})\right)^r &\leq &  \sum_{m=0}^\infty t^m \left(\sum_{k=m/2}^m
 x_k^r
 \exp(-\lambda r\kappa^{m-k})\right)
 \\
 &=&\sum_{k=0}^\infty x_k^r
 \left(\sum_{m=k}^{2k}t^m \exp(-\lambda r\kappa^{m-k})\right)\\
&=& \sum_{k=0}^\infty x_k^r t^k \left(\sum_{m=k}^{2k}
t^{m-k}\exp(-\lambda r\kappa^{m-k})\right)\\
& \leq&  \sum_{k=0}^\infty x_k^r t^k
\left(\sum_{m=0}^\infty t^m \exp(-\lambda r\kappa^m)\right)\\
&=& K(t,r,\kappa)\, \sum_{k=0}^\infty x_k^r t^r,
\end{eqnarray*}
because the series $\sum_{m=0}^\infty t^m \exp(-\lambda
r\kappa^{m})$ is convergent.

\noindent \textsc{Case 2. $r>1.$}\\ First, we extend the inner sum
to $k$ from $0$ to $m.$ Also, to avoid problems with summability,
we replace the infinite series $\sum_{m=0}^\infty$ with a finite
one $\sum_{m=0}^M,$ prove the appropriate inequality with a
constant $K$ not depending on $M$ and pass to the limit $M\to
\infty $ afterwards.

Let $a=\frac {2\ln t}{r\ln \kappa}$ (so that $\kappa^{ar}=t^2$).
There exists a constant $C>0$ such that $e^{-\lambda x}\leq
C\,x^{-a},$ for $x>0.$ It follows:
  \begin{eqnarray} \sum_{m=0}^M t^m \left(\sum_{k=0}^m x_k e^{-\lambda \kappa^{m-k}}
  \right)^r & \leq & C^r \sum_{m=0}^Mt^m \left(\sum_{k=0}^m x_k\,
  \frac{1}{\kappa^{(m-k)a}}\right)^r\nonumber\\
 &=& C^r \sum_{m=0}^M\frac{t^m}{\kappa^{amr}}
 \left(\sum_{k=0}^m x_k \kappa^{ka}\right)^r\nonumber\\
 &=& C^r \sum_{m=0}^M \frac{1}{t^m}\left(\sum_{k=0}^m
 x_k\,\kappa^{ka}\right)^r.\label{sss}
 \end{eqnarray}
Denote $S_m=\sum_{k=0}^m x_k \kappa^{ka}$ (for completeness, set
$S_{-1}=0$) and let $\tau=\frac{1}{t} (<1).$ We have:
\[
\sum_{m=0}^M \tau^m (S_m^r -S_{m-1}^r) =\sum_{m=1}^{M-1} S_m^r
\tau^m(1-\tau)+\tau^MS_M^r\geq (1-\tau)\sum_{m=0}^M \tau^mS_m^r,
\]
and so
\begin{eqnarray*} (\ref{sss}) &=& C^r\sum_{m=0}^M\tau^m S_m^r\\
&\leq& \frac{C^r}{1-\tau} \sum_{m=0}^M \tau^m
(S_m^r-S_{m-1}^r)\\
&\leq& \frac{rC^r}{1-\tau} \sum_{m=0}^M\tau^m
x_m\kappa^{ma}\,S_{m}^{r-1}
\end{eqnarray*}
(this is so because for $x,y >0$ we have $(x+y)^r-x^r\leq
ry(x+y)^{r-1}$).

\smallskip

Now use the  following discrete H\"{o}lder inequality:

\medskip

for $M=1,2,..., $ $\tau >0,$ $A_m, B_m \geq 0$ for
$m=0,1,2,...,M,$ and $p,q>1$ such that
$\frac{1}{p}+\frac{1}{q}=1,$
\begin{equation}\label{discrholder}
\sum_{m=0}^M A_m B_m\tau^m \leq \left(\sum_{m=0}^M A_m^p\tau^m
\right)^{\frac{1}{p}}\left(\sum_{m=0}^M
B_m^q\tau^m\right)^{\frac{1}{q}}.
\end{equation}
Applying (\ref{discrholder}) with $p=r,$ $q=\frac{r}{r-1},$
$A_m=x_m\kappa^{ma},$ $B_m=S_m^{r-1}$ we get:
\[(\ref{sss})\;\leq \frac{rC^r}{1-\tau}
 \left(\sum_{m=0}^M\tau^mx_m^r\kappa^{mar}\right)^{\frac{1}{r}}
\left(\sum_{m=0}^M \tau^m S_m^r\right)^{\frac{r-1}{r}}, \] which
results in
\[
\sum_{m=0}^M \tau^m S_m^r \leq \frac{rC^r}{1-\tau}
\left(\sum_{m=0}^M\tau^mx_m^r\kappa^{mar}\right)^{\frac{1}{r}}
\left(\sum_{m=0}^M \tau^m S_m^r\right)^{\frac{r-1}{r}},
\]
and further in
\[
\left(\sum_{m=0}^M \tau^m S_m^r\right)^\frac{1}{r} \leq
\frac{rC^r}{1-\tau}
\left(\sum_{m=0}^M\tau^mx_m^r\kappa^{mar}\right)^{\frac{1}{r}}.\]
Since $S_m=\sum_{k=0}^m x_k \kappa^{ka},$ $\tau=\frac{1}{t},$ and
$\kappa^{ar}=t^r,$ from (\ref{sss}) it follows that we are
done.\hfill $\Box$
\\

\noindent{\bf Conclusion of the proof of Theorem \ref{mmain}.} The
first sum in (\ref{ttt}) is estimated by (\ref{modhardy}), with
$t=2^{q(\alpha+\frac{d}{p})},$ $x_k=i_k^{(\alpha)}(f),$
$\kappa=2^{\frac{d_w}{d_w-1}}. $ For the second sum in
(\ref{ttt}), first forget about the exponential factor (which is
smaller than 1 anyway), and then use (\ref{clashardy}) with $t$
and $x_k$ as in the first sum. What we get is:
\[I_1^{(\alpha)}(f)\leq c\sum_{m=1}^\infty 2^{mq(\alpha+\frac{d}{p})}
(i_m^{(\alpha)}(f))^{\frac{q}{p}}=c\sum_{m=0}^\infty 2^{m\alpha q}
\left(2^{md}i_m^{(\alpha)}(f)\right)^{\frac{q}{p}}=\|(a_m^{(\alpha)}(f))\|_{q}^q.
\]

Collecting all the estimates obtained, we get that
\[I^{(\alpha)}(f)\leq c \,\|f\|_{L^p}^q + c\,\|(a_m(f))\|_q^q,\]
and further
\[(I^{(\alpha)}(f))^{1/q}\leq c\, \|f\|_{L^p} + c\, \|(a_m(f))\|_q.\]
This concludes the proof.\hfill $\Box$

This  theorem has a natural extension   to the case $q=\infty.$ In
For $\alpha=\frac{d_w}{2},$ $p=2$ and $q=\infty$ it has been
proven in \cite{Jon,KPP2,GHL}. To obtain the desired result, we
basically follow the lines of \cite{KPP2} and \cite{GHL}.

\begin{theo}\label{infty}
[extension to the case $q=\infty$] Let $\alpha>0, p\geq 1,$ $f\in
L^p(X,\mu)$ and let $a_m^{(\alpha)}(f)$  be defined as before.
Then \(\ds\sup_{m>0} a_m^{(\alpha)}(f)\) is finite if and only if
$$\sup_{t\in(0,1)}\frac{1}{t^{\frac{p\,\alpha}{d_w}}} \int_X\int_X
|f(x)-f(y)|^pp(t,x,y)\,d\mu(x)d\mu(y)<\infty.$$  Moreover,
\[\|(a_m^{(\alpha)}(f))\|_\infty \asymp \|f\|_p + \sup_{t>0}\frac{1}{t^{\frac{p\,\alpha}{d_w}}}
\int_X\int_X |f(x)-f(y)|^pp(t,x,y)\,d\mu(x)d\mu(y),\]
\end{theo}

Before we start the proof of Theorem \ref{infty}, let us state and
prove the following simple lemma.

\begin{lem}\label{suma}
Suppose $C,\alpha,\beta,\gamma >0.$ Then there exist constants
$K_1,K_2=K_{1,2}(C,\alpha,\beta,\gamma)$ such that
\begin{equation}\label{suma1}
K_1 \,t^{\alpha\beta}\leq
\sum_{m=0}^{\infty}2^{-m\alpha}\exp(-\left(\frac{C}{2^mt^\beta}\right)^\gamma)\leq
K_2\,t^{\alpha\beta}.
\end{equation}
\end{lem}

\noindent {\bf Proof of the lemma.} Consider
\begin{eqnarray*}
I_C&:=& \int_0^1
\exp(-\left(\frac{Cx^{\frac{1}{\alpha}}}{t^\beta}\right)^\gamma)\,dx
=\sum_{m=0}^\infty \int_{2^{-(m+1)\alpha}}^{2^{-m\alpha}}
\exp(-\left(\frac{Cx^{\frac{1}{\alpha}}}{t^\beta}\right)^\gamma)\,dx=:
\sum_{m=0}^\infty I_m.
\end{eqnarray*}

The function $x\mapsto
\exp(-\left(\frac{Cx^{\frac{1}{\alpha}}}{t^\beta}\right)^\gamma)$
is monotone decreasing, and so
\[\frac{1}{2^{m\alpha}}\left(1-\frac{1}{2^\alpha}\right)
\exp(-\left(\frac{C}{2^mt^\beta}\right)^\gamma) \leq I_m \leq
\frac{1}{2^{m\alpha}}\left(1-\frac{1}{2^\alpha}\right)
\exp(-\left(\frac{1}{2^\alpha}\frac{C}{2^mt^\beta}\right)^\gamma).
\]

Summing up over $m,$ we get
\[\left(1-\frac{1}{2^\alpha} \right) \sum_{m=0}^\infty
\frac{1}{2^{m\alpha}}\,
\exp(-\left(\frac{C}{2^mt^\beta}\right)^\gamma)\leq I_C\leq
\left(1-\frac{1}{2^\alpha} \right) \sum_{m=0}^\infty
\frac{1}{2^{m\alpha}}\,
\exp(-\left(\frac{1}{2^\alpha}\frac{C}{2^mt^\beta}\right)^\gamma)
\]
and it follows that
\[
\frac{2^\alpha}{2^\alpha-1}I_{C/2^\alpha} \leq
\sum_{m=0^\infty}2^{-m\alpha}
\exp(-\left(\frac{C}{2^mt^\beta}\right)^\gamma)\leq
\frac{2^\alpha}{2^\alpha-1}I_C.
\]
The last thing we have to do is to single out the dependence of
$I_C$ on $t.$ Substitute $z=x^{\gamma/\alpha}t^{-\beta\gamma}$ in
the integral, so that
\[I_C=t^{\alpha\beta}\int_0^{t^{-\beta\gamma}}\exp(-C^\gamma z)z^{\alpha/\gamma-1}dz,\]
which is integrable as long as $\alpha/\gamma >0.$ We are
done.\hfill $\Box$

\noindent{\bf Proof of Theorem \ref{infty}.} \\
 \noindent\textsc{Part 1.
The lower bound.} Suppose that $f\in L^p(X,\mu)$ is such that
$$s(f)=\sup_{t\in(0,1)}\frac{1}{t^{\frac{p\,\alpha}{d_w}}} \int_X\int_X
|f(x)-f(y)|^pp(t,x,y)\,d\mu(x)d\mu(y)<\infty.$$ Similarly to what
we have done before, restrict the area of integration to the set
$\{\rho(x,y)\leq t^{\frac{1}{d_w}}\},$ and use the bound for the
transition density, so that
\begin{eqnarray*}
s(f)& \geq&  c\,\sup_{t\in(0,1)}
\frac{1}{t^{\frac{p\,\alpha}{d_w}}}\frac{1}{t^{\frac{d}{d_w}}}
\int\int_{\rho(x,y)\leq
t^{\frac{1}{d_w}}}|f(x)-f(y)|^pd\mu(x)d\mu(y)\\
&\geq & c \sup_m 2^{mp\alpha} 2^{md} \int\int_{\rho(x,y)\leq
2^{-m}}|f(x)-f(y)|^pd\mu(x)d\mu(y),
\end{eqnarray*}
where the last inequality was obtained by using $t_m= 2^{-md_w}.$

\noindent \textsc{Part 2. The upper bound.}  Now suppose that
$f\in L^p$ and that $\|(a_m^{(\alpha)})(f)\|_\infty<\infty.$

Write the  integral in the definition of $s(f)$  as
 \begin{eqnarray*}
 &&\sum_{m=0}^\infty \int\int_{2^{-(m+1)}< \rho(x,y)\leq
2^{-m}} |f(x)-f(y)|^p\,p(t,x,y)\,d\mu(x) d\mu(y)\\
&\leq &c\, t^{-d/d_w}\sum_{m=0}^\infty \exp {(-c_4
(\frac{1}{2^mt^{1/d_w}})^\frac{d_w}{d_w-1})}
\int\int_{\rho(x,y)\leq 2^{-m}}|f(x)-f(y)|^p\,d\mu(x)d\mu(y)\\
&\leq & c\|(a_m^{(\alpha)})(f)\|_\infty
t^{-d/d_w}\sum_{m=0}^\infty 2^{-m(d+p\alpha)}\,\exp{(-c_4
(\frac{1}{2^mt^{1/d_w}})^{\frac{d_w}{d_w-1}})}.
\end{eqnarray*}

Lemma \ref{suma} allows us to estimate the last sum by $c\,
t^{\frac{d+p\alpha}{d_w}},$ and so  we get
\[ \frac{1}{t^{\frac{p\,\alpha}{d_w}}}\int_X\int_X
|f(x)-f(y)|^p\,p(t,x,y)\,d\mu(x)d\mu(y)  \leq
c\|(a_m^{(\alpha)}(f)\|_\infty.\label{ost}
\]
We are done.\hfill $\Box$

\begin{rem}[$q<\infty$]\label{rem1}
  {\rm For a function $f\in L^p(X,\mu)$ it is immediate to establish that
\begin{equation}\label{quaqua}
\int_1^\infty \frac{1}{t^{\frac{\alpha q}{d_w}}}
\left(\int_X\int_X |f(x)-f(y)|^p
p(t,x,y)\,d\mu(x)d\mu(y)\right)^{\frac{q}{p}}\frac{dt}{t}\leq
c\|f\|_{L^p}^q.
\end{equation}
Indeed, from symmetry we have
\begin{eqnarray*}
&&\int_X\int_X |f(x)-f(y)|^p p(t,x,y)\, d\mu(x)d\mu(y)\\
&\leq & 2^p \int_X \int_X |f(x)|^p p(t,x,y)\, d\mu(x)d\mu(y)\\
&=& 2^p \int_X |f(x)|^p (\int_X p(t,x,y)\,d\mu(y))d\mu(x)=2^p
\|f\|_{L^p}^p
\end{eqnarray*}
and (\ref{quaqua}) follows.

Therefore the norm in $Lip(\alpha,p,q)(X)$ is also equivalent to
\[\|f\|_{L_p}+(\widetilde{I}^{(\alpha)}(f))^{\frac{1}{q}},\]
where
$$\widetilde{I}^{(\alpha)}(f)= \int_0^\infty \frac{1}{t^{\frac{\alpha
q}{d_w}}}\left(\int_X\int_X |f(x)-f(y)|^p
p(t,x,y)\,d\mu(x)d\mu(y)\right)^{\frac{q}{p}} \frac{dt}{t}.$$

}
\end{rem}

\begin{rem}
 [$q=\infty$] {\rm Similarly, in this case we can take the supremum over $t>0$ instead of $t\in(0,1)$ and still
 get an equivalent norm.}
\end{rem}

\subsection{Range of parameters allowed}
The definition of Besov-Lipschitz spaces, as well as our
characterisation, work for arbitrary values of $\alpha>0,p,q\geq
1$. However, for some triples of parameters the resulting spaces
are trivial, and consist of constant functions only.

So far we have satisfactory results for $p=2$ only.
 It is known
that (see \cite{Jon}, Cor. 3 for the Sierpi\'{n}ski gasket,
\cite{KPP2}, Prop. 2 for the general case) the spaces $Lip(\alpha,
2, \infty)(X)$ are degenerate when $\alpha >\frac{d_w}{2}. $ From
here, it is immediate to see that $Lip(\alpha,2,q)(X)$ are
degenerate as well when $\alpha>\frac{d_w}{2}$. On the other hand,
when $\alpha\leq \frac{d_w}{2},$ then the spaces
$Lip(\alpha,2,\infty)(X)$ are dense in $L^2(X).$

Therefore $\alpha\leq \frac{d_w}{2}$ is the natural threshold for
the spaces $Lip(\alpha,2,\infty).$

Consider now the spaces $Lip(\alpha,p,p)(X).$ It is clear that for
$f\in L^p(X)$ we have
 \begin{equation}\label{lll}
 \|f\|_{\alpha,p,p}\asymp \int_X\int_X
\frac{|f(x)-f(y)|^p}{\rho(x,y)^{d+p\alpha}} \,d\mu(x)d\mu(y)
+\|f\|_{L^p}^p
\end{equation}
(when $\mbox{diam}\,X<\infty,$ then the term $\|f\|_{L^p}^p$ can
be omitted).

In \cite{KPP3} it has been proven that the finiteness of the
integral in (\ref{lll}), when $p=2$ and $\alpha\geq \frac{d_w}{2}$
implies that $f=const$ (and for $\alpha<\frac{d_w}{2}$ we get
dense subspaces of $L^2(X),$ which are domains of the {\em stable
processes} on $X,$ see \cite{Sto}). For the Sierpi\'{n}ski gasket,
this result (degeneracy) was earlier proved in \cite{Jon}. When
$\mbox{diam}\,X<\infty,$ then by an application of Jensen's
inequality we get the same conclusion for $p\geq 2$ (i.e.
$Lip(\alpha,p,p)$ degenerate when $\alpha>\frac{d_w}{2}$). We do
not know whether the value $\alpha=\frac{d_w}{2}$ is critical when
$p>2,$ and where should the threshold be placed when $1\leq p<2.$

For open subsets of the Euclidean space it is known that $\alpha=1
(=\frac{d_w(\mathbb{R}^d)}{2})$ works for all values of $p\geq 1,$
see \cite{Bre}. We do not expect this to hold in general.

\section{Links with other definitions of Besov-Lipschitz spaces}

\subsection{Strichartz Besov spaces on simple fractals}\label{stri}

Strichartz in \cite{Str} introduced the definition of three types
of Besov spaces on the Sierpi\'nski gasket first, and then on
p.c.f. self-similar fractals.

For the Sierpi\'nski gasket, the spaces
$(\Lambda^{p,q}_\alpha)^{(1)}(X)$ are defined for
$\frac{d}{p}<\alpha\leq \frac{d}{p}+\alpha_1,$ where
$\alpha_1=\frac{\log 2}{\log 5/3}=\frac{1}{d_w-d},$ and consist of
those bounded continuous functions on the gasket (this is why
there is a restriction to $\frac{d}{p}<\alpha$), for which
\begin{equation}\label{deff}
\|((2^{d_w-d})^{m\alpha}\delta_{m,p}(f))_m\|_{q}<\infty,
\end{equation}
 where
\[\delta_{m,p}(f)=\left\{
\begin{array}{ll}
 \left(2^{-md} \sum_{x\stackrel {m}{\sim} y} |f(x)-f(y)|^p\right)^\frac{1}{p}& \mbox{if } p<\infty,\\[2mm]
\sup\{|f(x)-f(y)|: {x\stackrel {m}{\sim} y} \} & \mbox{if }
p=\infty.
\end{array}
\right.
\]
The relation ${x\stackrel {m}{\sim} y} $ is the neighbouring
relation on the $m-$th approximation of the gasket, $V_m.$ Namely,
${x\stackrel {m}{\sim} y} $ if and only if $x,y\in V_m$ and
$|x-y|=\frac{1}{2^m}$ (for a precise definition we refer to
\cite{Str}).

Theorem 3.18 in \cite{Bod} and the discussion thereafter assert
that on simple fractals one has $(\Lambda^{p,q}_\alpha)^{(1)}(X)=
Lip (\alpha/(d_w-d), p,q)(X).$ Jonsson spaces do not carry the
restriction $\frac{d}{p}<\alpha,$ but on the other hand we know
that $Lip (\alpha,p,p)(X)$ are empty when $p\geq 2$ and
$\alpha\geq \frac{d_w}{2}$ (when $p>2,$ then the proof requires
that $\mbox{diam}\,X<\infty$). Therefore the admissible $\alpha$'s
in the  definition (\ref{deff}) should be restricted to
$\alpha\leq \frac{(d_w-d)d_w}{2}\approx 0.920042,$ which is a
smaller number than $\frac{d}{p}+\frac{1}{d_w-d}$ (e.g. when
$p=2,$ then $\frac{d}{2}+\frac{1}{d_w-d}\approx 2.14939665$).

Since Theorems \ref{mmain} and \ref{infty} give a characterisation
of the spaces $Lip(\alpha,p,q)(X),$ they simultaneously
characterise the spaces $(\Lambda^{p,q}_\alpha)^{(1)}(X).$

We do not whether the other spaces defined in \cite{Str}:
$(\Lambda^{p,q}_\alpha)^{(2)}(X)$ and
$(\Lambda^{p,q}_\alpha)^{(3)}(X)$  (defined through higher order
differences) can be characterised in a similar manner.

\subsection{Hu-Z\"{a}hle Besov spaces}\label{huza}

In \cite{HZ}, the authors introduce the following Besov-Lipschitz
spaces on $(X,\rho,\mu)$, under identical assumptions on $(X,\rho,\mu)$ as those  in the present paper:\\
for $\beta>0, p,q\geq 1,$
\[B_{\beta}^{p,q}(X)=\{f\in L^p(X,\mu): \left(\int_0^\infty(t^{k-\frac{\beta}{2}} \|\frac{\partial^k}{\partial t^k}P_t f\|_p)^q
\frac{dt}{t}\right)^{\frac{1}{q}}<\infty \},\] where
$k=\left[\frac{\beta}{2}\right]+1$ ($(P_t)_{t>0}$ is the semigroup
with kernel $p(\cdot,\cdot,\cdot)$). When $k=1$ (i.e. $\beta<2$),
they also establish that (see Theorem 5.2 of \cite{HZ})
$B_{\beta}^{2,2}(X)=H^{\beta}_{2}(X),$ with equivalent norms,
where
\[H_2^\beta(X)=\{ f\in L^2(X,\mu): \int_0^\infty (1+\lambda)^{{\beta}/{2}}d\langle E_\lambda
f,f\rangle<\infty\}.
\]
$(E_\lambda)_{\lambda>0}$ is the resolution of identity of the
generator of the semigroup $(P_t)_{t>0}$ of s.a. operators on
$L^2(X,\mu),$  (the actual definition was different, but this
equivalent condition is also given in \cite{HZ}).

We know that when $\beta<2,$ then $H^\beta_2(X)= Lip (\frac{\beta
d_w}{4},2,2),$ with equivalent norms. It follows from the earlier
results of St\'os \cite{Sto}, and also from the results in
\cite{HZ}. It can also be derived from our Theorem \ref{mmain}
together with Remark \ref{rem1}.

Indeed, the condition $f\in Lip(\frac{\beta d_w}{4},2,2)$ means
that
\begin{eqnarray*}\infty> \int_{0}^\infty
\frac{1}{t^{\beta/2}}\int_X\int_X |f(x)-f(y)|^2
p(t,x,y)\,\frac{dt}{t} &=& 2\int_0^\infty \frac{1}{t^{\beta/2}}
\,\frac{1}{t} \langle f-P_tf,f\rangle\, dt\\
&=& \int_0^\infty \int_0^\infty \frac{1-e^{-\lambda t}}{t}d\langle
E_\lambda f, f\rangle\,dt \\
&=& \int_0^\infty\left(\int_0^\infty \frac{1-e^{-\lambda
t}}{t^{1+\beta/2}}\,dt\right) d\langle E_\lambda f,f\rangle.
 \end{eqnarray*}

 After the substitution $s=\lambda t,$ the inner integral becomes
 $C_\beta\lambda^{\beta/2},$ with $C_\beta= \int_0^\infty
 \frac{1-e^{-s}}{s^{1+{\beta/2}}}\,ds,$ which is finite as long
 as $\beta <2.$
 And since we are dealing with $f\in L^2(X),$ the finiteness of
 $\int_0^\infty \lambda^{\beta/2}d\langle E_\lambda f,f
 \rangle$ implies the finiteness of $\int_0^\infty (1+\lambda)^{\beta/2}d\langle E_\lambda f,f
 \rangle,$ so that $f\in H^\beta_2(X).$

\medskip

\noindent {\bf Open question.} We do not know what is the relation
of spaces $Lip(\alpha,p,q)(X)$ to the spaces $B_\beta^{p,q}(X),$
in general -- even when $p=2.$ From  our characterisation we see
that $f\in Lip(\alpha,2,q)(X)$ if and only if $f\in L^2(X)$ and

\begin{eqnarray*}
\infty&\geq &
 \int_0^\infty \frac{1}{t^{\frac{\alpha
q}{d_w}}}\left(\int_X\int_X|f(x)-f(y)|^2p(t,x,y)\,d\mu(x)d\mu(y)\right)^{q/2}\,\frac{dt}{t}\\
&=& \int_0^\infty \frac{1}{t^{\frac{\alpha q}{d_w}}} \left(
\langle f-P_tf,f\rangle\right)^{\frac{q}{2}}\, \frac{dt}{t}
<\infty
\\
&=&
 \int_0^\infty \frac{1}{t^{\frac{\alpha
q}{d_w}-\frac{q}{2}}}\left(\frac{1}{t}\int_0^\infty
(1-e^{-t\lambda})d\langle E_\lambda
f,f\rangle\right)^{q/2}\frac{dt}{t}.
\end{eqnarray*}
On the other hand, $f\in B^\beta_{p,q}(X)$ if and only if
\[\int_0^\infty \left(t^{k-\frac{\beta}{2}}\left(\int_0^\infty \lambda^{2k}
 e^{-2\lambda t} d\langle E_\lambda f,f\rangle\right)^{1/2}\right)^{q}
\frac{dt}{t}<\infty\] (this is so because
$\|\frac{\partial^k}{\partial t^k}P_tf\|_2^2= \int_0^\infty
\lambda^{2k}e^{-2t\lambda}d\langle E_\lambda f,f\rangle$).

Guided by the case $q=2$ we expect to have some kind of
relationship between the two conditions when $\beta<2$ (i.e. when
$k=1$). It is unclear whether they yield the same sets of
functions, and what should be the dependence between $\beta$ and
$\alpha$ in general. When $p\neq 2,$ we can no longer use the
spectral representation, and so the situation is even more
unaccountable.

\end{document}